\documentclass[12pt,reqno,utf8]{amsart}

\usepackage{latexsym}
\usepackage{amsmath}
\usepackage{amssymb, amscd, amsthm, mathtools}
\usepackage[all]{xy}
\usepackage{multirow}
\usepackage{xcolor}
\usepackage{tikz}
\usepackage{tikz-cd}
\usepackage{mathrsfs}
\usepackage{fancyhdr}
\usepackage{cite}
\usepackage{hyperref}
\usepackage[a4paper, hcentering, left=3.5cm, right=3.5cm, top=3.5cm, bottom=3.5cm]{geometry}

\usepackage[shortlabels,inline]{enumitem}
\SetEnumerateShortLabel{a}{\textup{(\alph*)}}
\SetEnumerateShortLabel{A}{\textup{(\Alph*)}}
\SetEnumerateShortLabel{1}{\textup{(\arabic*)}}
\SetEnumerateShortLabel{i}{\textup{(\roman*)}}
\SetEnumerateShortLabel{I}{\textup{(\Roman*)}}

\DeclareMathOperator{\R}{\mathbb{R}}
\DeclareMathOperator{\bP}{\mathbb{P}}
\DeclareMathOperator{\A}{\mathbb{A}}

\DeclareMathOperator{\bO}{\mathcal{O}}

\DeclareMathOperator{\Conv}{\mathrm{Conv}}
\DeclareMathOperator{\Spec}{\mathrm{Spec}}

\DeclareMathOperator{\chara}{\mathrm{char}}

\DeclareMathOperator{\red}{\mathrm{red}}

\DeclareMathOperator{\Pic}{\mathrm{Pic}}
\DeclareMathOperator{\sm}{\mathrm{sm}}

\DeclareMathOperator{\Ker}{\mathrm{Ker}}

\DeclareMathOperator{\cl}{\mathrm{cl}}
\DeclareMathOperator{\anh}{\mathrm{h}}
\DeclareMathOperator{\inte}{\mathrm{int}}

\newtheorem{thm}{Theorem}[section]
\newtheorem{lem}[thm]{Lemma}
\newtheorem{prop}[thm]{Proposition}
\newtheorem{cor}[thm]{Corollary}
\newtheorem{conj}{Conjecture}

\theoremstyle{definition}

\theoremstyle{remark}
\newtheorem*{rmk}{Remark}

\begin{document}
\title[Bertini-type theorems over real closed fields]{Generalized connectedness and Bertini-type theorems over real closed fields}
\author{Yi Ouyang$^{1,2}$ and Chenhao Zhang$^1$}
\address{$^1$School of Mathematical Sciences, Wu Wen-Tsun Key Laboratory of Mathematics,   University of Science and Technology of China, Hefei 230026, Anhui, China}

\address{$^2$Hefei National Laboratory, University of Science and Technology of China, Hefei 230088, China}

\email{yiouyang@ustc.edu.cn}	
\email{chhzh@mail.ustc.edu.cn}


\thanks{Partially supported by NSFC (Grant No. 12371013) and  Innovation Program for Quantum Science and Technology (Grant No. 2021ZD0302902).}

\begin{abstract}
    In this paper, we establish a real closed analogue of Bertini's theorem. 
    Let $R$ be a real closed field and $X$ a formally real integral algebraic variety over $R$. We show that if the zero locus of a nonzero global section $s$ of an invertible sheaf on $X$  has a formally real generic point, then $s$ does not change sign on $X$, and vice versa under certain conditions. As a consequence,  we demonstrate that there exists a nonempty open subset of hypersurface sections preserving formal reality and integrality for quasi-projective varieties of dimension $\geq 2$ under these conditions. 
\end{abstract}
\maketitle
% Introduction
\section{Introduction} For a smooth projective variety $Y\subseteq \bP^n_k$, the classical Bertini theorem states that a general hyperplane $H \subseteq \Gamma(\bP^n_k, \bO(1))$ intersects $Y$ in a smooth subscheme (see \cite[II, Theorem~8.18]{Hart77}) if $k$ is an algebraically closed field. If $k$ is a finite field, Poonen established the existence of a hypersurface $H$ in $\bP^n_k$ such that $H \cap Y$ is smooth (see \cite{Poon04}).

Let $R$ be a real closed field and $X$ a formally real integral algebraic variety over $R$. In this paper
we develop an analogue of Bertini's theorem over $R$. 
Suppose $s$ is a nonzero global section of an invertible sheaf $\mathcal{L}$ on $X$. Our main theorem states that $s$ does not change sign on $X$ if  its zero locus $V(s)$  has a formally real generic point, and vice versa  under regularity assumptions and assuming a certain conjecture (Conjecture~\ref{Conj_M}) holds for curves. 
Based on this result,  we derive  the following Bertini-type results: 
\begin{enumerate}
    \item Let $G=\{0\neq s\in\Gamma(X,\mathcal{L})\mid V(s)\ \text{has a formally real generic point}\}$. If Conjecture~\ref{Conj_M} holds for curves over $R$, then for any vector subspace $L \subseteq \Gamma(X,\mathcal{L})$ of finite dimension $\geq 2$, $G \cap L$ has nonempty interior in $L$ (under the order topology).
        This holds unconditionally if $R$ is replaced by an archimedean (but not necessarily real closed) field.
    \item If Conjecture~\ref{Conj_M} holds for curves over $R$, then there exists a nonempty open subset of hypersurface sections preserving formal reality and integrality for a quasi-projective variety of dimension $\geq 2$ over $R$.
\end{enumerate}

\section{Preliminaries}
In this paper, $R$ is always a real closed field.
It is known that $\chara(R)=0$ and $R$ admits a unique ordering compatible with its field structure (refer to \cite{Lam05}). 
Let $R^+=\{x\in R\mid x>0\}$.

For an algebraic variety $X$ over $R$, let $X^{\anh}$ be the topological space defined on the set $X(R)$ induced by the order topology of $R$ in \cite[Proposition~3.1]{Bria12}, whose dimension $\dim(X^{\anh})$ is defined to be the Krull dimension of $X(R)$.
For a morphism $f \colon X \to Y$  of $R$-algebraic varieties, we denote by $f^{\anh}: X^{\anh}\to Y^{\anh}$ its induced continuous morphism.

A subset $T \subseteq R^n$ is called convex if for all $p,q \in T$ and $\lambda \in [0,1]$, we have $\lambda p + (1-\lambda)q \in T$.
Let $S_1 \subseteq R^{n_1}$ and $S_2 \subseteq R^{n_2}$.
A map $\phi \colon S_1 \to S_2$ is called a convex map if $\phi$ maps convex subsets of $S_1$ to convex subsets of $S_2$.
We denote by $\Conv(S_1,S_2)$ the set of all continuous convex maps from $S_1$ to $S_2$.

\begin{lem}[\cite{Van98}, Chapter~1, Lemma~3.6] \label{Lem_BCoAF}
    Let $f \in R(x_1,\dots,x_n)$.
    If $f$ is regular on a convex subset $T$ of $R^n$, then $f$ maps $T$ to a convex subset of $R$.
\end{lem}

\begin{lem}[\cite{Van98}, Chapter~7, \S2] \label{Lem_MVToRF}
    Let $f \in R(x)$ and $f'$ be its derivative  with respect to $x$.
    Let $T$ be a convex subset of $R$ in which $f$ has no singularities. Then 
    \begin{enumerate}[1]
      \item $f \in R$ if and only if $f'$ takes the value $0$ at infinitely many points.
      \item If $f \notin R$, then $f'|_T\geq 0$ if and only if $f$ is increasing on $T$.
      \item If $f \notin R$, then $f'|_T\leq 0$ if and only if $f$ is decreasing on $T$.
      \item If $a < b \in T$ and $f(a)=f(b)$, then $f'$ has a zero in $[a,b]$.
    \end{enumerate}
\end{lem}

Let $S$ be a subset of $R$ and let $a \in R$.
We write $a < S$ if $a < s$ for every $s \in S$.
The notations $a > S$, $a \leq S$, and $a \geq S$ are defined analogously.
\begin{lem} \label{Lem_PoIFT}
    Let $\phi \colon X \to Y$ be a morphism of algebraic varieties over $R$.
    If $\phi$ is \'{e}tale at $p \in X(R)$, then $\phi^{\anh}$ is a local homeomorphism at $p \in X^{\anh}$.
\end{lem}
\begin{proof}
    The problem is local.
    Without loss of generality, assume $X=\Spec(B)$, $Y=\Spec(A)$, and $\phi$ is standard \'{e}tale.
    Let $A=R[x_1,\dots,x_m]/I$, $B=A[x]_Q/P$.
    We can view $Y^{\anh}$ as a closed subspace of $R^m$, and $X^{\anh}$ as a subspace of $R^{m+1}$.
    We regard $P$ and $Q$ as polynomials in $R[x_1,\dots,x_m,x]$.
    Since $\phi$ is smooth at $p$, we have $P'_x(p) \neq 0$.
    Since $x \mapsto -x$ is an automorphism of $R[x]$, we may assume $P'_x(p) > 0$.
    Let $V(Q)$ be the zero set of $Q$ in $R^{m+1}$, which is a closed subset of $R^{m+1}$.
    Then by the continuity of $P'_x$ and Lemma~\ref{Lem_BCoAF}, there exists an open convex subset $T_0$ of $R^m$ and an open interval $S_0$ of $R$ such that $p \in T_0 \times S_0$, $T_0 \times S_0 \cap V(Q)=\emptyset$ and $0 < P'_x(T_0 \times S_0)$.
    
    Let the coordinates of $p$ be $(\phi(p),p_0)$.
    By Lemma~\ref{Lem_MVToRF}, $P(\phi(p),x)$ is increasing on $S_0$.
    Take $a_1<a_2 \in S_0$ such that $P(\phi(p),a_1)<0$ and $P(\phi(p),a_2)>0$.
    By the definition of $S_0$, we certainly have $p \in T_0\times (a_1,a_2)$.
    By continuity, there exists a convex open neighborhood $T_1 \subseteq T_0$ of $\phi(p)$ such that $P(T_1,a_1)<0$ and $P(T_1,a_2)>0$.
    By Lemma~\ref{Lem_BCoAF} and Lemma~\ref{Lem_MVToRF}(2), for all $y \in T_1$, $\left(\left.\phi^{\anh}\right|_{T_1\times (a_1,a_2) \cap X^{\anh}}\right)^{-1}(y)$ are singletons.
    Therefore, \[\phi^{\anh} \colon T_1\times (a_1,a_2) \cap X^{\anh} \to T_1 \cap Y^{\anh}\] is a bijection.
    Since the projection map $ T_1\times (a_1,a_2) \to T_1$ is an open map, it follows that $\left.\phi^{\anh}\right|_{T_1\times (a_1,a_2) \cap X^{\anh}}$ is also an open map.
    Hence, $\phi^{\anh}$ is a homeomorphism between $T_1\times (a_1,a_2) \cap X^{\anh}$ and $T_1 \cap Y^{\anh}$.
\end{proof}

\begin{lem} \label{Lem_OToAS}
    For $f \in R[x_1,\dots,x_n] \setminus\{0\}$, $V(f)^{\anh}$ is nowhere dense in $R^n$.
\end{lem}
\begin{proof}
    If $n=0$, it is obvious. Assume $n > 0$.
    Suppose there exists a nonempty open set \[S=(a_1,b_1) \times (a_2,b_2) \times \dots \times (a_n,b_n)\] such that $f(S)=\{0\}$, where $(a_i,b_i)$ are open intervals in $R$.
    Let \[S_0=(a_1,b_1) \times \dots \times (a_{n-1},b_{n-1}).\]
    
    Let $f=\sum_{j=1}^m f_jx_n^j$, where $f_j \in R[x_1,\dots,x_{n-1}]$.
    For all $q \in S_0$, $f(q,(a_n,b_n))=\{0\}$.
    Since $(a_n,b_n)$ contains infinitely many elements, $f(q,x_n)=0$.
    Therefore, for all $j=1,\dots,m$, $f_j(S_0)=\{0\}$.
    By induction on dimension, we conclude that $f_j=0$.
\end{proof}

\begin{prop} \label{Prop_DToAS} For
  an algebraic variety $X$ over $R$,  $\dim(X^{\anh})\geq n$ if and only if there exist an open subscheme $U$ of $X$ and a morphism $\phi \colon U \to \A^n_R$ of $R$-algebraic varieties such that $\phi^{\anh}(U^{\anh})$ is not nowhere dense in $R^n$.
\end{prop}
\begin{proof}
    Let $Z$ be the Zariski closure of the $R$-points in $X$.
    
\medskip \noindent Proof of $\Rightarrow$.
    By \cite[\href{https://stacks.math.columbia.edu/tag/00OT}{Tag 00OT}]{stacks-project}, we may assume $X$ is an affine integral scheme.
    Let $X=\Spec(A)$, $Z=\Spec(B)$, and $A \twoheadrightarrow B$ be the surjection corresponding to the closed immersion $Z \hookrightarrow X$.
    Suppose there exists a principal open subscheme $V_0=\Spec(B_g)$ of $Z$ and a morphism $\phi_0 \colon V_0 \to \A^n_R$ of $R$-algebraic varieties such that $\phi_0^{\anh}(V_0^{\anh})$ is not nowhere dense in $R^n$, where $g$ is an element in $A$.
    By the universal property of polynomial rings, the morphism $R[x_1,\dots,x_n] \to B_{g}$ can be lifted to $A_{g}$, and the corresponding scheme morphism $\phi \colon V \to \A^n_R$ satisfies the requirement.
    Therefore, we may assume $Z=X$.
     
    Choose an irreducible component $X_0$ of $X$ with maximal dimension. 
    By \cite[Theorem 4.1.4]{Boch98}, the generic points of $X$ are formally real.
    By \cite[Theorem 4.1.2]{Boch98}, $\dim(X_0^{\anh})\geq n$.
    Therefore, we may assume $X$ is a smooth integral scheme.
    
    Let $p \in X(R)$.
    By \cite[\href{https://stacks.math.columbia.edu/tag/054L}{Tag 054L}]{stacks-project}, there exist a Zariski open neighborhood $U$ of $p$ and $\Phi \colon U \to \A^{\dim(X)}_R$ such that $\Phi$ is \'{e}tale at $p$.
    By Lemma~\ref{Lem_PoIFT}, $\Phi^{\anh}$ is a local homeomorphism at $p$.
    Let $\beta \colon \A^{\dim(X)} \twoheadrightarrow \A^n$ be a projection morphism.
    We define $\phi=\beta\circ \Phi$.
    Since both the projection morphism and a homeomorphism are open maps, it follows that $\phi^{\anh}$ is an open map at the point $p$.
    
\medskip \noindent Proof of  $\Leftarrow$.
    Let $W$ be the Zariski closure of $\phi^{\anh}(U^{\anh})$ in $\A^n$.
    By \cite[\href{https://stacks.math.columbia.edu/tag/00P1}{Tag 00P1}]{stacks-project}, we have \[\dim(W) \leq \dim(U\cap Z)\leq \dim(Z).\]
    Since $\phi^{\anh}(U^{\anh})$ is not nowhere dense in $R^n$, $W=\A^n$ by Lemma~\ref{Lem_OToAS}.
    Therefore, $\dim(X^{\anh}) \geq n$.
\end{proof}

\section{Main results}
Let $Y$ be an $R$-algebraic variety and $V \subseteq Y^{\anh}$.
We say $V$ is a generalized connected subset of $Y^{\anh}$ if:
\begin{enumerate}[i]
    \item There exist an open affine subset $U \subseteq Y^{\red}$ containing $V$ and a smooth morphism $\phi \colon U \to \A^n_R$;
    \item There exists an open set $V_0 \subseteq U^{\anh}$ containing $V$ such that $\phi^{\anh}|_{V_0}$ is an open embedding;
    \item $\phi^{\anh}(V)$ is convex in $R^n$, and for every $f \in \Gamma(V,\bO_X)$, the map $f\circ (\phi^{\anh}|_{V})^{-1} \colon \phi^{\anh}(V) \to R$ is convex.
\end{enumerate}
\begin{conj} \label{Conj_M}
    For a smooth algebraic variety $X$ over $R$, the space $X^{\anh}$ always admits a covering by generalized connected open subsets.
\end{conj}

\begin{prop} \label{Prop_CoFRF}
 Suppose $R$ is a real closed field. 
    \begin{enumerate}[1]
      \item Conjecture~\ref{Conj_M} holds for curves over $R$ if $\Conv(R^+,R)$ is a subring of $C(R^+,R)$.
      \item Conjecture~\ref{Conj_M} holds in general if $R$ is moreover archimedean.
    \end{enumerate}
\end{prop}
\begin{proof}
    (1). Let $X$ be a smooth algebraic curve over $R$.
    Let $p\in X^h$, we need to find a generalized connected open set containing $p$.
    By the definition of generalized connected sets, we may assume $X=\Spec(R[x,y]_h/g)$ and $g$ is irreducible.
    Since $X$ is smooth, we have $(g'_x(p),g'_y(p))\neq 0$.
    Using affine transformations on $k^2$, we may assume $g'_x(p) \neq 0$ and $g'_y(p) \neq 0$.
    
    There exists a convex open neighborhood $T_0 \subseteq k^2$ of $p$ such that
    $0 \notin g'_x(T_0)$, $0 \notin g'_y(T_0)$ and \[T_0 \cap \Spec(R[x,y]/(g,h)) =\emptyset.\]
    Let the coordinates of $p$ in $k^2$ be $(p_0,p_1)$.
    Using the same method as in the proof of Lemma~\ref{Lem_PoIFT}, we can find open intervals $(a_1,a_2)$ and $(b_1,b_2)$ satisfying:
    \begin{enumerate}[i]
      \item $p \in (a_1,a_2)\times (b_1,b_2) \subseteq T_0$;
      \item $g$ maintains constant sign on $(a_1,a_2) \times\{b_1\}$ and $(a_1,a_2) \times\{b_2\}$.
      \item $g$ maintains constant sign on $\{a_1\} \times(b_1,b_2)$ and $\{a_2\} \times(b_1,b_2)$.
    \end{enumerate}
    Let $P_x$ be the projection of the curve $X$ onto the $x$-axis, and let $P_y$ be the projection onto the $y$-axis.
    Then $P_x|_{(a_1,a_2)\times(b_1,b_2)}$ and $P_y|_{(a_1,a_2)\times(b_1,b_2)}$ are both open embeddings.
    Let $\phi=P_y \circ (P_x)^{-1}|_{(a_1,a_2)}$.
    Then $\phi$ is a monotonic continuous convex map from $(a_1,a_2)$ to $R$.
    
    Let $f=\frac{f_1}{f_2} \in \Gamma((a_1,a_2)\times(b_1,b_2), \bO_X)$, where $f_1,f_2 \in R[x,y]$ are coprime.
    Since $\Spec(R[x,y] / (f_1,f_2))$ is finite, there exists a finite open cover $\{(\gamma_{1j},\gamma_{2j})\}_{j=1}^m$ of the interval $(a_1,a_2)$
    such that on each interval $(\gamma_{1j},\gamma_{2j})$, $f=\frac{f_{1j}}{f_{2j}}$ and $f_{2j} \in R[x,y]$ has no zeros in $(\gamma_{1j},\gamma_{2j})\times(b_1,b_2) \cap X^{\anh}$.
    The union of intersecting convex subsets remains convex, and the gluing of continuous convex maps on convex subsets remains convex.
    Therefore, we only need to prove that for all $j=1,\dots,m$, $f\circ (P_x)^{-1}|_{(\gamma_{1j},\gamma_{2j})}$ is a continuous convex map.
    
    The nonempty open interval $(\gamma_1,\gamma_2)$ is homeomorphic to $(0,+\infty)$ via the continuous convex map $\frac{1}{x-\gamma_1}-\frac{1}{\gamma_2-\gamma_1}$.
    Thus, $\Conv((\gamma_1,\gamma_2),R)$ forms a ring.
    Therefore, $f_{1j}\circ (P_x)^{-1}|_{(\gamma_{1j},\gamma_{2j})}$ and $f_{2j}\circ (P_x)^{-1}|_{(\gamma_{1j},\gamma_{2j})}$ are both convex maps.
    The inverse of a convex subset of $R$ not containing $0$ remains a convex subset of $R$, so 
    $\frac{1}{f_{2j}}\circ (P_x)^{-1}|_{(\gamma_{1j},\gamma_{2j})}$ is a convex map.
    In conclusion, $f\circ (P_x)^{-1}|_{(\gamma_{1j},\gamma_{2j})}$ is a continuous convex map.
    
    (2). Connected subsets of $\R$ are convex; thus, Conjecture~\ref{Conj_M} holds for $\R$.
    Since $R$ is archimedean, $R$ is isomorphic to a subfield of $\R$.
    Let $X$ be a smooth algebraic variety over $R$, and let $p \in X^{\anh}$.
    Then there exist a Zariski open neighborhood $U$ of $p$, an \'{e}tale morphism $\phi \colon U \to \A^m_R$, and a neighborhood $T \subseteq U^{\anh}$ such that $\phi^{\anh}$ restricts to an open embedding on $T$ with $S:=\phi^{\anh}(T)$ being a convex subset of $R$.
    The map $\phi^h$ being algebraic automatically extends to $U(\R)$.
    For any $f \in \Gamma(T)$, $f\circ(\phi^{\anh}|_T)^{-1}$ is a continuous map on $\inte(\cl_{\R^m}(S))$.
    Therefore, $f\circ(\phi^{\anh}|_T)^{-1}|_{S}$ is a continuous convex map.
\end{proof}
\begin{rmk}
 It is known that $\Conv(\R^+,\R)=C(\R^+,\R)$. 
 Thus $\Conv(\R^+,\R)$ forms a ring. It would be interesting to ask the following questions: Is $\Conv(R^+,R)$ always a subring of $C(R^+,R)$? If not, under which conditions is $\Conv(R^+,R)$ a subring?
\end{rmk}

\begin{thm} \label{Thm_CoSaR}
 Let $X$ be a formally real integral algebraic variety over $R$.  
 \begin{enumerate}[1]
    	\item For $\mathcal{L} \in \Pic(X)$ and $0 \neq s \in \Gamma(X, \mathcal{L})$, if no generic point of $V(s)$ is formally real, then $s$ does not change sign on $X$. Namely,
    	there exists an affine trivialization of $\mathcal{L}$ restricted to $X_{\sm}$
    	\[\{U_i,\left.\mathcal{L}\right| _{U_i}\xrightarrow{\sim}\mathcal{O}_{U_i}\}_{i\in I},\]
    	such that either $s(T)\geq 0$ or $s(T)\leq 0$ for every generalized connected subset $T$ of $(U_i)^{\anh}$.
    	\item Assume Conjecture~\ref{Conj_M} holds for curves over $R$. Then the converse of (1) is true if $X$ is regular in codimension one and $V(s)$ is a reduced subscheme of $X$.
    \end{enumerate}
\end{thm}
\begin{proof}
    (1). If $s$ changes sign on $X$, then there exists an affine open subscheme $U \subseteq X_{\sm}$ and a trivialization $\left.\mathcal{L}\right|_{U} \simeq \bO_U$ such that $s$ changes sign on a generalized connected subset $T_0 \subseteq U^{\anh}$.
    
    Let $n=\dim(X)$, and let $\phi \colon U_0 \to \A^n_R$ be the morphism satisfying the requirements in the definition of generalized connected sets, where $U_0$ is an affine open subset of $U$.
    We first prove that there exists $p_0\in T_0$ such that $s$ changes sign on every open neighborhood of $p_0$.
    Since $s$ changes sign on $T_0$, $T_0$ is not a singleton.
    Take $p_1,p_2 \in T_0$ such that $s(p_1)<0$ and $s(p_2)>0$.
    Consider the affine line $\A_R^1$ passing through $\phi(p_1)$ and $\phi(p_2)$.
    Let $a_1$ be the coordinate of $p_1$ in $\A_R^1$, and $a_2$ the coordinate of $p_2$; without loss of generality, assume $a_1<a_2$.
    Let $Y$ be the pullback of $U_0$ along this affine line.
    Denote $(\phi^{\anh}|_{T_0})^{-1}([a_1,a_2])$ by $T_1$.
    By definition, $s$ changes sign on $T_1$.
    
    Let $\psi=\left.(\phi^{\anh}|_{T_0})^{-1}\right|_{[a_1,a_2]} \colon [a_1,a_2] \to T_1$.
    Define
    \[S_1=\{a\in [a_1,a_2]\mid \exists a_0 \in [a_1,a_2],\ s\circ\psi(a_0)\leq 0\  \text{and}\ a_0\geq a\}.\]
    Since $s$ is not identically zero on $Y$, $s\circ\psi$ has only finitely many zeros.
    Therefore, there exists a maximal zero $z_0$ of $s\circ \psi$ in $S_1$.
    Let $S_2=(z_0,+\infty) \cap S_1$.
    Since $s\circ \psi$ is continuous, $S_1 \subsetneqq [a_1,a_2)$.
    Let $S_3=[a_1,a_2] \setminus S_1$.
    If $z_0$ is the supremum of $S_1$, then the claim is proved.
    If not, by continuity, $S_2$ has no supremum.
    But by definition, $S_2 \cup S_3$ is a convex subset of $[a_1,a_2]$, so $0 \in s\circ \psi(S_2 \cup S_3)$, a contradiction.
    Let $p_0 = (\phi^{\anh}|_{T_0})^{-1}(z_0) \in T_0$.
    In conclusion, $s$ changes sign in every open neighborhood of $p_0$.
    
    Without loss of generality, assume \[U_0=\Spec(R[x_1,\dots,x_{n+1}]_{h}/g).\]
    Let $B=R[x_1,\dots,x_{n+1}]_{h}/g$, $A=B/s$.
    The kernel $\Ker(R[x_1,\dots, x_n] \to A)$ is a principal ideal of $R[x_1,\dots, x_n]$, denoted by $(f)$.
    Therefore, $f=s\cdot s_1 \mod g$, where $s_1 \in R[x_1,\dots,x_{n+1}]$.
    Since $\phi$ is smooth, $s_1(p_0) \neq 0$.
    Hence, $s_1$ maintains constant sign near $p_0$ by continuity.
    Therefore, $f$ changes sign in every open neighborhood of $p_0$ (by taking the first $n$ coordinates of $p_0$).
    
    In summary, there exist $q_1$ and $q_2$ in $\phi(T_0)$ such that $f(q_1)=b_1 <0$ and $f(q_2)=b_2 >0$.
    Let $W$ be the hyperplane in $R^n$ perpendicular to the line through $q_1$ and $q_2$.
    By definition, there is an open subset $G_0 \subseteq U_0^{\anh}$ containing $T_0$, such that $\phi|_{G_0}$ is an open embedding and $\phi(G_0)$ is convex.
    Then by continuity, there exists a nonempty convex neighborhood $W_0\subseteq W$ of $0$ such that $f(q_1+W_0)<0$, $f(q_2+W_0)>0$, $q_1+W_0 \subseteq \phi(G_0)$, and $q_2+W_0 \subseteq \phi(G_0)$.
    
    Let $Z=\Spec(R[x_1,\dots, x_n]/f)$.
    Using an affine transformation to adjust coordinates, project $Z$ from the direction of the line through $q_1,q_2$ onto $W$.
    Denote this projection map by $\phi_1\colon Z \to \A^{n-1}_R$.
    By Lemma~\ref{Lem_BCoAF}, the image of $\phi_1^{\anh}$ contains $W_0$.
    Therefore, the image of the composite morphism $\Spec(A)^{\anh} \to Z^{\anh} \to R^{n-1}$ contains $W_0$.
    By Proposition~\ref{Prop_DToAS}, $\dim(V(s)^{\anh})\geq n-1$.
    Since $\dim(V(s))=n-1$, $V(s)$ has a formally real generic point.
    
    (2). If $V(s)$ is reduced and $V(s)$ has a formally real generic point, then there exists a formally real $R$-point $p$ in $V(s)_{\sm}$ .
    Since $X$ is regular in codimension one, we have $p \in X_{\sm}$.
    Let $U$ be an affine open neighborhood of $p$ in $X_{\sm}$, such that $\mathcal{L}$ has a trivialization $\left.\mathcal{L}\right|_{U} \simeq \bO_U$ and $\phi \colon U \to \A^n$ is an \'{e}tale morphism, where $n=\dim(X)$.
    
    We may assume $U=\Spec(R[x_1,\dots,x_{n+1}]_{h}/g)$.
    Let $B=R[x_1,\dots,x_{n+1}]_{h}/g$, $A=B/s$.
    $\Ker(R[x_1,\dots, x_n] \to A)$ is the principal ideal $(f)$ of $R[x_1,\dots, x_n]$.
    Let $C=R[x_1,\dots, x_n]/f$.
    Let $\mathfrak{m}$ be the maximal ideal in $R[x_1,\dots, x_n]$ corresponding to $\phi(p)$ in $\A^n_R$.
    Let $\mathfrak{n}$ be the maximal ideal in $B$ corresponding to $p$ in $\Spec(B)$.
    We have an injective homomorphism of local rings $\psi \colon C_{\mathfrak{m}} \to A_{\mathfrak{n}}$, where $A_{\mathfrak{n}}$ is a regular local ring.
    Let \[A'=C_{\mathfrak{m}} \otimes_{R[x_1,\dots, x_n]} B_{\mathfrak{n}}.\]
    Then $C_{\mathfrak{m}} \to A'$ is an \'{e}tale homomorphism.
    
    Note that $A_{\mathfrak{n}}$ is a quotient of $A'$; let $I =\Ker(A' \twoheadrightarrow A_{\mathfrak{n}})$.
    Since $A$ is reduced, $C$ is also reduced, so $A'$ is reduced.
    Suppose $\Spec(A')$ has only one irreducible component; then $A'=A_{\mathfrak{n}}$ is regular and $C_{\mathfrak{m}}$ is regular.
    By \cite[\href{https://stacks.math.columbia.edu/tag/00OF}{Tag 00OF}]{stacks-project},
    $C_{\mathfrak{m}}$ is regular if and only if the vector $(f'_{x_1}(p),\dots,f'_{x_n}(p))$ is nonzero.
    Without loss of generality, assume $f'_{x_1}(p) \neq 0$.
    Let the coordinates of $p$ be $(a_1,\dots,a_n)$.
    Then there exists $(\beta_1,\beta_2) \ni a_1$ such that $f'_{x_1}$ maintains constant sign on
    \[(\beta_1,\beta_2) \times \{a_2\} \times \dots \times \{a_n\}.\]
    Let $Y$ be the pullback of $\phi$ along the morphism $\A_R^1 \hookrightarrow \A_R^n$ (this morphism is the spectrum version of the quotient homomorphism $R[x_1,\dots,x_n] \twoheadrightarrow R[x_1]$).
    By assumption, there exists a generalized connected open subset $T_0$ of $Y$ such that $p \in T_0$ and $\phi|_Y(T_0) \subseteq (\beta_1,\beta_2)$.
    By the assumption on the interval $(\beta_1,\beta_2)$, $s$ changes sign on $T_0$.
    
    We now prove that $\Spec(A')$ has only one irreducible component, i.e., $I=0$.
    Since $C_{\mathfrak{m}} \to A'$ is \'{e}tale and $p$ is an $R$-point, we have
    \[R=C_{\mathfrak{m}}/\mathfrak{m}\xrightarrow{\sim}A'/\mathfrak{m}\xrightarrow{\sim}A_{\mathfrak{n}}/\mathfrak{n}.\]
    Therefore, $\mathfrak{m}A'=\mathfrak{n}A'$ and $\mathfrak{m}A_{\mathfrak{n}}=\mathfrak{n}A_{\mathfrak{n}}$.
    Let $M$ be the quotient $A_{\mathfrak{n}}/C_{\mathfrak{m}}$ as a $C_{\mathfrak{m}}$-module.
    Similarly, let $M'=A'/C_{\mathfrak{m}}$.
    We have $\mathfrak{m}M=0$ and $\mathfrak{m}M'=0$.
    Let $M_0=\Ker(M' \to M)$.
    We have a commutative diagram of exact sequences:
    \[\begin{tikzcd} &0 \arrow[r] \arrow[d] &I \arrow[r] \arrow[d] &M_0\arrow[d] & \\
        0 \arrow[r] &C_{\mathfrak{m}} \arrow[r]  \arrow[d]&A' \arrow[r]  \arrow[d]&M' \arrow[r]  \arrow[d]&0 \\
        0 \arrow[r] &C_{\mathfrak{m}} \arrow[r]  \arrow[d]&A_{\mathfrak{n}} \arrow[r]  \arrow[d]&M \arrow[r]  \arrow[d]&0 \\
         & 0 \arrow[r] & 0 \arrow[r] & 0  & \end{tikzcd}\]
    By the snake lemma, we have $I \xrightarrow{\sim} M_0$.
    Therefore, $\mathfrak{n}I=\mathfrak{m}I=0$.
    By Nakayama's lemma, $I=0$.
\end{proof}

\begin{prop} \label{Prop_FRoBT}
    Let $X$ be a formally real integral algebraic variety over $R$. 
    For $\mathcal{L}\in \Pic(X)$, let \[ G=\{0\neq s\in\Gamma(X,\mathcal{L})\mid V(s)\ \text{has a formally real generic point}\}.\]
    If Conjecture~\ref{Conj_M} holds for curves over $R$, then
    for any vector subspace $L \subseteq \Gamma(X,\mathcal{L})$ of finite dimension $\geq 2$, $G \cap L$ has nonempty interior in $L$ (under the order topology).
\end{prop}
\begin{proof}
    Let $U \subseteq X_{\sm}$ be a nonempty affine open subscheme such that $\mathcal{L}$ has a trivialization $\left.\mathcal{L}\right|_{U} \simeq \mathcal{O}_U$ on $U$.
    Let $\{s_i\}_{i=1,\dots,n}$ be an $R$-basis of $L$.
    Let $U'=U_{s_1s_2}$.
    Since $X$ is integral, $s_1\neq 0$ and $s_2 \neq 0$, we have $U'\neq \emptyset$.
    Since $X$ is formally real, there exists $p \in U'(R)$.
    Let $U_1=U'_{s_1(p)s_2-s_2(p)s_1}$.
    Similarly, $U_1 \neq \emptyset$.
    Since Conjecture~\ref{Conj_M} holds for curves, the union of generalized connected subsets containing $p$ is Zariski dense in $U'$.
    Therefore, there exists a generalized connected subset $T_0$ of $(U')^{\anh}$ containing $p$, such that $T_0\cap U_1^{\anh} \neq \emptyset$.
    Let $q\in T_0\cap U_1^{\anh}$.
    
    In summary, the matrix \[\begin{pmatrix}s_1(p) & s_2(p) \\s_1(q) & s_2(q) \end{pmatrix}\] is invertible.
    Therefore, there exist $k_1, k_2 \in R$ such that $(k_1s_1+k_2s_2)(p)<0$ and $(k_1s_1+k_2s_2)(q)>0$.
    
    The functions $e_p(x)=\sum_{i=1}^{n}s_i(p)x_i$ and $e_q(x)=\sum_{i=1}^{n}s_i(q)x_i$ are linear, hence continuous on $R^n$. For $v=(k_1,k_2,0,\dots,0)$, we have $e_p(v) < 0$ and $e_q(v) > 0$. 
    Therefore, there exists an open subset $W \subseteq L$ containing $v$ such that $e_p(W)\subseteq (-\infty, 0)$ and $e_q(W)\subseteq (0, +\infty)$.
    By Proposition~\ref{Thm_CoSaR}, for every $s \in W$, $V(s)$ has a formally real generic point.
\end{proof}
\begin{cor} \label{Cor_AFoBT} The result in Proposition~\ref{Prop_FRoBT} holds unconditionally if
 the field $R$  is replaced by an archimedean  ( but not necessarily real closed) field $k$.
\end{cor}
\begin{proof}
    Conjecture~\ref{Conj_M} holds for $\R$.
    Since $k$ is dense in $\R$, nonempty open subsets of $L\otimes_{k}\R$ contain nonempty open subsets of $L$.
\end{proof}

\begin{cor}\label{Cor_BToSFR}   
    Let $X \hookrightarrow \bP^n_{R}$ be a formally real integral quasi-projective variety over $R$ of dimension $\geq 2$.
    Let
    \[ \begin{split}
    	&W=\{0\neq H\in \Gamma(\bP^n_{R}, \bO_{\bP^n}(d))\mid X\cap H\ \text{is formally real and integral}\}, \\ 
    	&G=\{0\neq s\in \Gamma(X, \bO_X(d))\mid V(s)\ \text{is formally real and integral}\}.
    \end{split}\] 
If Conjecture~\ref{Conj_M} holds for curves over $R$, then
\begin{enumerate}[1]
	\item  $W$ has nonempty interior in $\Gamma(\bP^n_{R}, \bO_{\bP^n}(d))$.
 \item	For any vector subspace $L \subseteq \Gamma(X, \bO_X(d))$ of finite dimension $\geq 3$, $G \cap L$ has nonempty interior in $L$.
\end{enumerate}
\end{cor}
\begin{proof}
    Let $C$ be the algebraic closure of $R$. 
    By \cite[Chapter~VIII, Theorem~2.5]{Lam05}, $[C:R]=2$.
    Let $K$ be the function field of $X$.
    If $K \otimes_{R} C = K_1 \times K_2$ is reducible, then $K_1=K_2=K$.
    Then $C \supseteq K$, so $K$ is not formally real, a contradiction.
    By \cite[\href{https://stacks.math.columbia.edu/tag/054Q}{Tag 054Q}]{stacks-project}, $X$ is geometrically irreducible.
    
    Hypersurfaces that have smooth proper intersection with $X_{\mathrm{sm}}$ are general in $\Gamma(\bP^n_{R}, \bO_{\bP^n}(d))$ (see \cite[\href{https://stacks.math.columbia.edu/tag/0FD6}{Tag 0FD6}]{stacks-project} or \cite[II, Theorem~8.18]{Hart77}).
    We can choose a hypersurface $H_1 \in \Gamma(\bP^n_{R}, \bO_{\bP^n}(d))$ such that $X_1=X_{\sm} \cap H_1$ is a smooth variety of dimension $\dim(X)-1$.
    Since $\dim(X)\geq 2$, we can choose $H_2 \in \Gamma(\bP^n_{R}, \bO_{\bP^n}(d))$ such that $X_2=X_1\cap H_2$ is a smooth variety of dimension $\dim(X)-2$.
    If $\dim(X)>2$, we can choose $H_3 \in \Gamma(\bP^n_{R}, \bO_{\bP^n}(d))$ such that $X_2\cap H_3$ is a smooth variety of dimension $\dim(X)-3$.
    If $\dim(X)=2$, we can choose $H_3 \in \Gamma(\bP^n_{R}, \bO_{\bP^n}(d))$ such that $X_2\cap H_3= \emptyset$.
    It can be shown that $H_1,H_2,H_3$ are linearly independent in $\Gamma(X, \bO_X(d))$.
    
    In summary, hypersurfaces that have geometrically irreducible intersection with $X$ are general in $\Gamma(\bP^n_{R}, \bO_{\bP^n}(d))$ (see \cite[\href{https://stacks.math.columbia.edu/tag/0G4F}{Tag 0G4F}]{stacks-project}).
    By Lemma~\ref{Lem_OToAS}, Zariski open subsets of a vector space are dense in the order topology.
    Therefore, by Proposition~\ref{Prop_FRoBT}, $W$ has nonempty interior.
    The proof for $\Gamma(X, \bO_X(d))$ is the same as the above argument.
\end{proof}

\section*{Statements and Declarations}
\subsection*{Conflict of interest}
The authors declare that they have no conflict of interest.
\subsection*{Data Accessibility}
No datasets were generated or analyzed during this study. Data sharing is not applicable to this purely theoretical work.

\end{document}